\documentclass[reqno]{amsart}
\usepackage{hyperref}

\numberwithin{equation}{section}

\usepackage{amsmath, amsthm, amssymb, color}

\usepackage[T2A]{fontenc}

\usepackage[cp1251]{inputenc}

\usepackage[english]{babel}
\usepackage{url}
\usepackage[normalem]{ulem} 

\newtheorem{defn}{Definition}[section]

\newtheorem{thm}[defn]{Theorem}

\theoremstyle{remark}

\newtheorem{ex}[defn]{Example}

\newcommand{\h}{\mathcal{H}}
\newcommand{\mn}{\mathbb N}

\def\<{\langle}
\def\>{\rangle}
\newcommand{\la}{\langle}
\newcommand{\ra}{\rangle}

\newcommand{\seq[1]}{(#1_n)_{n\in I}}

\def\newin {\,\kern-0.4em\in\kern-0.15em}

\newcommand{\seqN[1]}{(#1_n)_{n\in \mathbb N}}

\usepackage{makecell} 

\begin{document}
\title[Constructions of dual frames compensating for erasures]
{Constructions of dual frames compensating for erasures with implementation}

\author[Lj. Aramba\v si\' c, D. T. Stoeva\hfil \hfilneg] {Ljiljana Aramba\v si\' c, Diana T.  Stoeva}

\address{Ljiljana Aramba\v si\' c \newline Department of Mathematics, Faculty of Science, University of Zagreb, Bijeni\v cka cesta 30, 10000 Zagreb, Croatia}
\email{arambas@math.hr}

\address{Diana T. Stoeva \newline
Faculty of Mathematics, University of Vienna,
Oskar-Morgenstern-Platz 1, 1090 Vienna, Austria}
\email{diana.stoeva@univie.ac.at}

 \keywords{frame, erasure, reconstruction, dual frame, canonical dual \\ \, \\
\indent This is a preprint of the following chapter: Ljiljana Aramba\v si\' c, Diana T. Stoeva, {\it Constructions of dual frames compensating for erasures}, accepted for publication in ``Women in Analysis and PDE'' (\url{https://link.springer.com/book/9783031570049}), edited by M. Chatzakou, M. Ruzhansky, and D. Stoeva, part of the sub-series  ``Research Perspectives Ghent Analysis and PDE Center'' of the book series ``Trends in Mahematics'', Birkh\"auser Cham, expected to appear in June 2024. The preprint is posted on arXiv with permission of the Publisher.}

\dedicatory{In memory of the second author's brother, Dimo Stoev, who was a wonderful person and a great specialist in his work} 

\maketitle

\begin{abstract} 
Let $I\subseteq \Bbb N$ be a finite or infinite set and let $\seq[x]$ be a frame for a separable Hilbert space $\h$. Consider transmission of a signal $h\in\h$ where a finite subset $(\<h,x_n\>)_{n\in E}$ of the frame coefficients $(\<h,x_n\>)_{n\in I}$ is lost. 
There are several approaches in the literature aiming recovery of $h$.
In this paper we focus on the approach based on construction of a dual frame of the reduced frame $(x_n)_{n\in I\setminus E}$ which is then used for perfect reconstruction from the preserved frame coefficients $(\<h,x_n\>)_{n\in I\setminus E}$.
There are several methods for such construction, starting from the canonical dual or any other dual frame of $\seq[x]$.
We implemented the algorithms for these methods and performed tests to compare their computational efficiency.
\end{abstract}


\section{Introduction and Notation} \label{intr}

Throughout the paper, $\h$ denotes a separable (finite or infinite-dimensional) Hilbert space. 
Let $I\subseteq \Bbb N$ be a finite or infinite set. A sequence $\seq[x]$ in $\h$ 
 is called a \emph{frame for $\h$} \cite{DS} if there exist positive constants $A$ and $B$ such that
\begin{equation}\label{frame def}
A\|x\|^2\leq \sum_{n\in I}|\la x,x_n\ra |^2\leq B\|x\|^2,\quad x \in \h.
\end{equation}
Let us recall some basic facts about frames needed for the paper.
To a given a frame $\seq[x]$ for $\h$, one associates the so called \emph{analysis operator} $U: \h\rightarrow \ell^2(I)$ determined by $Ux=(\langle x,x_n\rangle)_{n\in I},\,x\in \h$. The adjoint operator $U^*$ (called the \emph{synthesis operator}) can be expressed by $U^*((c_n)_{n\in I})=\sum_{n\in I}c_nx_n$ for $(c_n)_{n\in I}\in \ell^2(I)$. The composition operator $U^*U$ is called the \emph{frame operator} and it is invertible on $\h$ (i.e., it is bounded and bijective from $\h$ onto $\h$).

Given a frame $(x_n)_{n\in I}$ for $\h$, there always exists a frame $\seq[z]$ for $\h$ so that the reconstruction formula
\begin{equation*}\label{defaltdual}
x=\sum_{n\in I}\la x,x_n\ra z_n,\quad x \in \h,
\end{equation*}
holds; such $\seq[z]$ is called a \emph{dual frame of $\seq[x]$}.
One dual frame of $\seq[x]$ is given by $((U^*U)^{-1}x_n)_{n\in I}$; it is called the \emph{canonical dual of $\seq[x]$}.
When the frame $(x_n)_{n\in I}$ is at the same time a Schauder basis for $\h$, the canonical dual is the only dual frame of $(x_n)_{n\in I}$,
while a frame $(x_n)_{n\in I}$ that is not a Schauder basis for $\h$ (called an \emph{overcomplete} or \emph{redundant frame})
possesses infinitely many dual frames.
For more on general frame theory we refer e.g. to \cite{CKu, Cbook, Hbook, KC}.

The possibility for redundancy is very important for reconstruction purposes in applications (in particular, in signal and image processing, data compression, and other areas). Transmitting a signal $h\in\h$  (more precisely, transmitting the frame coefficients $(\<h,x_n\>)_{n\in I}$) via some redundant frame $\seq[x]$ for $\h$, allows perfect reconstruction of $h$ even if some of the frame coefficients gets lost or damaged during the transmission, assuming that the preserved frame coefficients  arise from frame elements which span the space $\h.$
 (If a non-redundant frame is used in this setting, 
 then loss of a frame coefficient would prevent from perfect recovery of $h$.)
Related to this, we consider the following notion introduced in \cite{LS}.
It is said that a finite set of indices $E\subsetneq I$ satisfies the \emph{minimal redundancy condition} (in short, the \emph{MRC}) \emph{for} a frame $(x_n)_{n\in I}$ for $\h$ if the linear span of the set $\{x_n:n\in E^c(=I\setminus E)\}$ is dense in $\h,$ that is, if
\begin{equation}\label{def_mrp}
\overline{\text{span}}\,\{x_n:n\in E^c\}=\h.
\end{equation}
As it is observed in \cite{LS}, based on \cite[Theorem~5.4.7]{Cbook},
in the above setting the expression \eqref{def_mrp} holds if and only if $(x_n)_{n\in E^c}$ is a frame for $\h$.
Thus, in the case of such a set $E$,
one can recover $h$ through the formula
\begin{equation}\label{approach1}
h=\sum_{n\in E^c}\la h,x_n\ra v_n
\end{equation}
via a dual frame $(v_n)_{n\in E^c}$ of $(x_n)_{n\in E^c}$.
This approach for recovery of $h$ was considered in \cite{AB,AD}, where several ways of construction of a dual frame of $(x_n)_{n\in E^c}$ were studied (in \cite{AB} the authors constructed the canonical dual of $(x_n)_{n\in E^c}$ starting with the canonical dual of $(x_n)_{n\in I}$; motivated by
\cite[Proposition~2.8 and Remark~2.9]{AB} and \cite{sc}, in \cite{AD} the authors used an arbitrary dual frame of $(x_n)_{n\in I}$ and obtained a dual frame of $(x_n)_{n\in E^c}$ that is not necessarily the canonical one).
Note that one can compute the canonical dual of $(x_n)_{n\in E^c}$ using the standard definition via the inverse of the respective frame operator, but this might not be very efficient computationally in high dimensional spaces as it involves inversion of $r\times r$ matrix where $r$ is the dimension of the space. That is why more efficient methods are of interest - e.g.
Theorem \ref{thm2.5}(ii), where only $k\times k$ matrix is inverted ($k$ being the cardinality of the erasure set and expected to be much smaller than $r$); Theorem \ref{thm2.5}(iii) which is based on iterative procedure and thus also expected to be more efficient than inverting an $r\times r$-matrix.

For other approaches for recovery of $h$ we refer to \cite{CK, HLSS, HanSun, HP, LS, LH, LoH}.

\vspace{.1in}
In this paper we focus on some of the constructions in \cite{AB,AD}, implementing the algorithms and examining 
 their computational efficiency.
The paper is organized as follows.

In Section~\ref{section2} we consider the canonical case, that is, we start with a frame $\seq[x]$ and its canonical dual $\seq[y]$.
Based on the results obtained in \cite{AB} (see also \cite{AD}), we summarize three ways for construction of the canonical dual of the reduced frame $(x_n)_{n\in E^c}$ after erasure via a set $E$ satisfying the MRC. Further, in the case of a finite frame $\seq[x]$, we implemented two of these ways, the iterative algorithm (Theorem \ref{thm2.5}(iii)) and the algorithm with matrix inversion (Theorem \ref{thm2.5}(ii)), and tested them for efficiency in time-computation,  more precisely, we compared the time for computation of the canonical dual of $(x_n)_{n\in E^c}$ via the codes based on these two algorithms and via the MATLAB \emph{pinv}-function that can also serve the purpose for determining the canonical dual. Results from the tests are given in Section~\ref{sec4}.

Section~\ref{section3} is devoted to the non-canonical case - we start with an arbitrary dual frame $\seq[z]$ for $\seq[x]$.
We use results from \cite{AD} which show that, with some additional assumptions, we can adapt formulas obtained in the canonical case to get a (not necessarily canonical) dual frame for $(x_n)_{n\in E^c}$. Again, we implemented  the iterative algorithm (Theorem \ref{thm2.5_arb}(i))  and the algorithm with matrix inversion (Theorem \ref{thm2.5_arb}(ii)) and tested them for efficiency in time-computation. Results from the tests are given in Section~\ref{sec4}. We also compare the time-efficiency of the algorithms from the canonical and the non-canonical case (Table 1 in Section~\ref{sec4}).

\vspace{.1in}
We end the section with some notation.
The identity operator on $\h$ and the corresponding identity matrix in the finite dimensional case  are denoted by $I$.
For $x,y\in \h$, $\theta_{y,x}$ denotes
 the rank one operator defined by $\theta_{y,x}(h)=\la h,x\ra y,$ which is  clearly a bounded operator on $\h$.

\section{Constructions based on the canonical dual}\label{section2}

Let $(x_n)_{n\in I}$ be a frame for $\h$ and let $E$ denote the indices of the erased frame coefficients. 
Here we consider ways for construction of the canonical dual of the reduced frame $(x_n)_{n\in E^c}$ starting with the canonical dual of $(x_n)_{n\in I}$.
In the following theorem we summarize some results obtained in
\cite[Proposition~2.4, Theorems~2.5, 2.12 and 2.14, and (2.29)]{AB}, see also \cite[Proposition 2.2]{AD}. 

\begin{thm} \label{thm2.5}
Let $(x_n)_{n\in I}$ be a frame for $\h$ and $(y_n)_{n\in I}$ its canonical dual frame.
Suppose that a finite set of indices $E=\{1,2,\ldots,k\}\subsetneq I$ satisfies the MRC for $(x_n)_{n\in I}$. Then the following holds:
\begin{enumerate}

 \item[{\rm (i)}] The canonical dual $(v_n)_{n\in E^C}$ of $(x_n)_{n\in E^C}$  can be written as
\begin{equation}\label{candual-operator}
v_n=(I-\sum_{i=1}^k\theta_{y_{i},x_{i}})^{-1}y_n,\quad n\in E^c.
\end{equation}

\item[{\rm (ii)}]
For each $n \in E^c$,   the numbers $\alpha_{n1}, \alpha_{n2},\ldots, \alpha_{nk},$
given  by the formula
\begin{equation}\label{XVIII}
{\small
\left[
\begin{array}{c}
\alpha_{n1}\\
\alpha_{n2}\\
\vdots\\
\alpha_{nk}
\end{array}
\right]=
\left(\left[\begin{array}{cccc}
\langle y_{1},x_{1}\rangle&\langle y_{2},x_{1}\rangle&\ldots&\langle y_{k},x_{1}\rangle\\
\langle y_{1},x_{2}\rangle&\langle y_{2},x_{2}\rangle&\ldots&\langle y_{k},x_{2}\rangle\\
\vdots&\vdots& &\vdots\\
\langle y_{1},x_{k}\rangle&\langle y_{2},x_{k}\rangle&\ldots&\langle y_{k},x_{k}\rangle
\end{array}
\right]-I\right)^{-1}
\left[
\begin{array}{c}
\langle y_n,x_{1}\rangle\\
\langle y_n,x_{2}\rangle\\
\vdots\\
\langle y_n,x_{k}\rangle
\end{array}
\right]
}
\end{equation}
are well defined and
the sequence $(v_n)_{n\in E^C}$ determined by
\begin{equation}\label{candual-matrix}
v_n=y_n-\sum_{i=1}^k\alpha_{ni}y_{i},\quad n\in E^c,
\end{equation}
  is the canonical dual of $(x_n)_{n\in E^C}$.

\item[{\rm (iii)}] Denote
$$v_n^0:=y_n, \quad n\in I.$$
Let  $E_j=\{1,2,\ldots,j\}$ for $j=1,2,\ldots,k.$   For $j$ from $1$ to $k$ let
\begin{eqnarray}
\label{candual-iterations-alpha_1} \alpha_n^j:&=&\frac{\la v_n^{j-1},x_j\ra}{1-\la v_j^{j-1},x_j\ra}, \quad n\in E_j^c,\\
\label{candual-iterations_1}v_n^j:&=&v_n^{j-1}+ \alpha_n^j v_j^{j-1}, \quad 
n\in E_j^c.
\end{eqnarray}
Then, for every $j$ from $1$ to $k$, the sequences in \eqref{candual-iterations-alpha_1} and \eqref{candual-iterations_1} are well defined, and the sequence $(v_n^j)_{n\in  E_j^c}$ is the canonical dual of $(x_n)_{n\in E_j^c}.$
\end{enumerate}

\end{thm}

The part (iii) from the preceding theorem gives us an iterative procedure for defining the canonical dual of $(x_n)_{n\in E^c}.$
Moreover,  for each $j=1,\ldots,k,$  the $j$-th iteration gives the canonical dual frame of $(x_n)_{n\in E_j^c}.$

\section{Constructions based on an arbitrary dual frame}\label{section3}

In \cite{AD} the authors investigated whether Theorem~\ref{thm2.5} holds if one replaces the canonical dual $\seq[y]$ with an arbitrary dual $\seq[z]$ of a given frame $\seq[x].$ It turns out that in this more general case it can happen that some of the formulas from the canonical case are not well defined.
 Let $\seq[x]$ be a frame for $\h$ and let $\seq[z]$ be an arbitrary dual frame of $\seq[x]$. Consider a finite set of indices  $E=\{1,2,\ldots,k\}{\subsetneq I}$ that satisfies the MRC for $\seq[x]$ and denote
\begin{equation}\label{matrix-xz}
A_{X,Z,E}=\left[\begin{array}{cccc}
\langle z_{1},x_{1}\rangle&\langle z_{2},x_{1}\rangle&\ldots&\langle z_{k},x_{1}\rangle\\
\langle z_{1},x_{2}\rangle&\langle z_{2},x_{2}\rangle&\ldots&\langle z_{k},x_{2}\rangle\\
\vdots&\vdots& &\vdots\\
\langle z_{1},x_{k}\rangle&\langle z_{2},x_{k}\rangle&\ldots&\langle z_{k},x_{k}\rangle
\end{array}
\right]-I.
\end{equation}
In the case  when $\seq[z]$ is the canonical dual of $\seq[x]$, the matrix $A_{X,Z,E}$ appears in Theorem~\ref{thm2.5}(ii) and it is necessarily invertible.
The next example from \cite{AD} shows that there are (noncanonical) dual frames for which the associated matrix \eqref{matrix-xz} is not invertible. The same example also serves to show that the operator
$I-\sum_{i=1}^k\theta_{z_{i},x_{i}} $ from \eqref{candual-operator} need not be invertible, and that the adapted iterative procedure from Theorem~\ref{thm2.5}(iii) need not be well defined.

\begin{ex}\label{abc_not_d} 
Let $\seqN[e]$ be an orthonormal basis of $\h.$ Consider the frame
$\seqN[x]=(e_1, e_1, e_1, e_2, e_3, e_4, \ldots)$ and its (non-canonical) dual frame
$$ \seqN[z]=(e_1, -\frac{1}{2} e_1, \frac{1}{2}e_1, e_2, e_3, e_4,\ldots).$$
Obviously, the set $E=\{1\}$ satisfies the MRC for  $\seqN[x]$. However:

(a) The matrix $A_{X,Z,E}$ is the zero matrix, so it is not invertible.

(b) The operator $I-\theta_{z_1,x_1}=I-\theta_{e_1,e_1}$ is not invertible, because $e_1$ belongs to its kernel.

(c) If we take $v_n^0=z_n,$ $n\in \mn$, then the first step in \eqref{candual-iterations-alpha} cannot be done since  $\la z_1,x_1\ra=1$, so $\alpha_n^1=\frac{\la z_n,x_1\ra}{1-\la z_1,x_1\ra}$ for $n\geq 2$ is not well defined.
\end{ex}

In the following theorem we give a non-canonical version of Theorem~\ref{thm2.5}.

\begin{thm} \label{thm2.5_arb}
Let $\seq[x]$ be a frame for $\h$ and $\seq[z]$ an arbitrary dual frame  of $\seq[x]$.
Suppose that a finite set of indices $E=\{1,2,\ldots,k\}\subsetneq I$ satisfies the MRC for $(x_n)_{n\in I}$. Let  $E_j=\{1,2,\ldots,j\}$ for $j=1,2,\ldots,k.$
Denote
$$v_n^0:=z_n, \quad n\in I.$$
Suppose that, for $j$ from $1$ to $k$, the following sequences are well defined:
\begin{eqnarray}
\label{candual-iterations-alpha} \alpha_n^j:&=&\frac{\la v_n^{j-1},x_j\ra}{1-\la v_j^{j-1},x_j\ra}, \quad n\in E_j^c,\\
\label{candual-iterations}v_n^j:&=&v_n^{j-1}+ \alpha_n^j v_j^{j-1}, \quad n\in E_j^c.
\end{eqnarray}
Then the following statements hold:
\begin{enumerate}
\item[{\rm (i)}]
For every $j$ from $1$ to $k$, the sequence $(v_n^j)_{n\in  E_j^c}$ is a dual frame of $(x_n)_{n\in E_j^c}.$

\item[{\rm (ii)}]
The matrix $A_{X,Z,E}$ is invertible and the sequence
$(v_n)_{n\in E^C}$, determined by
\begin{equation}\label{candual-matrix-arb}
v_n=z_n-\sum_{i=1}^k\alpha_{ni}z_{i},\quad n\in E^c,
\end{equation}
where, for each $n \in E^c$,   the numbers $\alpha_{n1}, \alpha_{n2},\ldots, \alpha_{nk}$ are
given  by
\begin{equation}\label{XVIII-arb}
{\small
\left[
\alpha_{n1}\
\alpha_{n2}\
\ldots\
\alpha_{nk}
\right]^T= A_{X,Z,E}^{-1}
\left[
\langle z_n,x_{1}\rangle\
\langle z_n,x_{2}\rangle\
\ldots\
\langle z_n,x_{k}\rangle
\right]^T
,}
\end{equation}
is a dual frame of $(x_n)_{n\in E^C}$.

 \item[{\rm (iii)}] The operator $I-\sum_{i=1}^k\theta_{z_i,x_i}$ is invertible and the sequence $(v_n)_{n\in E^C}$ determined by
\begin{equation}\label{candual-operator-arb}
v_n=(I-\sum_{i=1}^k\theta_{z_{i},x_{i}})^{-1}z_n,\quad n\in E^c,
\end{equation}
is a dual frame of $(x_n)_{n\in E^C}$.

\item[{\rm (iv)}] The dual frames $(v_n)_{n\in E^C}$ of $(x_n)_{n\in E^C}$,
constructed in {\rm(i)}, {\rm (ii)}, and {\rm (iii)}, are the same. 
\end{enumerate}
\end{thm}

Note that, while validity of the iterative procedure determined by (\ref{candual-iterations-alpha})-(\ref{candual-iterations}) implies validity of the other two constructions (in (ii) and (iii)), there is no equivalence in this respect - there exist cases where (ii) and (iii) apply, but the iterative procedure does not apply, see \cite{AD} for more details.

\section{Implementation and computational efficiency}\label{sec4}

In this section we examine
the computational efficiency of the approaches in Sections~\ref{section2}  and \ref{section3}.
In the finite-dimensional case, we have implemented the algorithms of  Theorem~\ref{thm2.5_arb}(i)-(ii) (resp. Theorem~\ref{thm2.5} (iii)-(ii)), which provide constructions of a dual frame (resp. the canonical dual) of the reduced frame $(x_n)_{n\in E^C}$.
The scripts are available on \sloppy \url{http://dtstoeva.podserver.info/DualFramesCompensatingErasures.html}. 
The programming is done under the MATLAB environment, using also frame-commands from LTAFT\footnote{The Large Time-Frequency Analysis Toolbox (a Matlab/Octave open source toolbox for dealing with time-frequency analysis and synthesis),
\url{http://ltfat.org/}, see e.g. \cite{ltfatnote030,ltfatnote015}.}.

We tested the efficiency of the scripts for various frames $X=(x_n)_{n=1}^N$ varying the number $N$ of the frame elements, the dimension $r$ of the space, the redundancy $N/r$,
 the cardinality $k$ of the erasure set $E=\{1,2,\ldots,k\}$, as well as the starting dual frame $Z$ of $X$.
The elapsed time recorded in Table 1 is in seconds, measured
 using the MATLAB tic-toc functions.
The tests were performed on
Lenovo Thinkpad T450, 
Processor Intel(R) Core(TM) i5-5200U CPU @ 2.20GHz   2.19 GHz, 64-bit operating system, RAM 8GB, under Windows 10 Pro.

We compare the time for computing the canonical dual of the reduced frame $(x_n)_{n\in E^c}$ via the code
based on Theorem~\ref{thm2.5}(iii) ($t_1$ in Table 1),
via the code based on  Theorem \ref{thm2.5}(ii) ($t_2$  in Table 1), and via the pseudo-inverse\footnote{The pseudo-inverse approach is based on the fact that the synthesis operator of the canonical dual of a frame $X$ is the adjoint of the Moore-Penrose pseudoinverse of the synthesis operator of $X$ \cite[Theorem 1.6.6]{Cbook}.}
approach  ($t_3$ in Table 1)  using the MATLAB \emph{pinv}-function\footnote{Our first aim was to do comparison with the LTFAT function \emph{framedual} for computing the canonical dual, but since for general frames \emph{framedual} uses the pseudo-inverse approach calling the \emph{pinv}-function from MATLAB, we compare directly to the \emph{pinv}-function to avoid unnecessary delay.}.  
For the same frame $X$, for which the aforementioned tests were performed concerning the canonical dual,
we also compare the time  for computing another dual frame of the reduced frame $(x_n)_{n\in E^c}$,  
starting with an arbitrary dual frame $Z_i$ of $X$,  and using the iterative method in Theorem~\ref{thm2.5_arb}(i)
($t_4(Z_i)$ in Table 1) and  the 
 formulas (\ref{candual-matrix-arb})-(\ref{XVIII-arb})
 in Theorem~\ref{thm2.5_arb}(ii) ($t_5(Z_i)$ in Table 1), $i=1,2$.

In Table~1 we present some results from the tests - the execution time of the considered procedures and the respective errors. For each procedure in a test, the respective
error $e_i$
is computed using the MATLAB 2-norm function of
the matrix $V^* U-I$, where $U$ denotes the analysis operator of  $(x_n)_{n\in E^c}$ and $V$ denotes  the analysis operator of the constructed dual  frame $(v_n)_{n\in E^c}$ via the considered procedure.
Through the results of the tests, on the one hand, one can compare the performance of the two procedures
 -  via the iterative algorithm and via inversion of the matrix $A_{X,Z,E}$ (as well as comparison to the third procedure in the canonical case - the use of the MATLAB \emph{pinv}-function),
and on the other hand, one can compare the performance of the canonical dual versus  another dual frame for the desired constructions. 
In each test, the shortest executed time is marked in blue, and for each of the three sub-groups $t_1-t_3$, $t_4(Z_1)-t_5(Z_1)$, and $t_4(Z_2)-t_5(Z_2)$ of the respective test, the shortest time in the sub-group is marked with bold style.\
        The scripts used to produce the tests reflected in Table 1 
         are available on
 \url{http://dtstoeva.podserver.info/DualFramesCompensatingErasures.html}.

Concerning Tests 1-8: With $N$ and $r$  determined by the user as input  parameters, 
the test-program produces a frame $X=(x_n)_{n=1}^N$ with random elements, the canonical dual of $X$, and 
  two randomly determined dual frames $Z_1$ and $Z_2$ of $X$. The user can enter a value of $k$ as wished until the program verifies that $E$ has the MRC for $X$. For this set $E$ the program measures $t_1$-$t_3$, 
 $t_4(Z_1)$ and $t_5(Z_1)$, $t_4(Z_2)$ and $t_5(Z_2)$.

Concerning Test 9: 
With $N$ and $r$ initialized at the beginning of the script 
(the values can be easily changed by 
 the user), 
a specific simple frame 
$X=(x_n)_{n=1}^N$ and its canonical dual are constructed, a specific non-canonical dual frame $Z_1$ of $X$ is also constructed (aiming to include a case where Theorem~\ref{thm2.5_arb}(i)-(ii) does not apply), and another dual frame $Z_2$ of $X$ is randomly chosen (aiming to include   
a non-canonical case where Theorem~\ref{thm2.5_arb}(i)-(ii)  applies).  
 Running the program, the user can enter a value of $k$ as wished until the program verifies that $E$ has the MRC for $X$. 
For this set $E$ the program measures $t_1$-$t_3$, $t_4(Z_1)$-$t_5(Z_1)$ (for the chosen value $k=4$ the iterative procedure of Theorem~\ref{thm2.5_arb}(i) for $Z_1$ cannot be completed and Theorem~\ref{thm2.5_arb}(ii) does not apply), and $t_4(Z_2)$-$t_5(Z_2)$.

\begin{table}[h!] \label{tab:table1new}
  \begin{center}
        {\tiny
          \begin{tabular}{|l|r|r|r|r|r|r|r|r|r|}
         \hline
       \ & Test 1  & Test 2 & Test 3 & Test 4 & Test 5 & Test 6 & Test 7 & Test 8 & Test 9\\
          \hline
     \hspace{-.09in}  $N$\hspace{-.09in}  & 6000 & 6000 & 6000 & 7000 & 5000 & 8000 & 8000 & 8000 &
      3010 \\
               \hline
     \hspace{-.09in}  $r$ \hspace{-.09in} & 4000 &4000& 4000 & 4000 & 4000 & 200  & 2000 & 6000 & 3000 \\
               \hline
     \hspace{-.09in}   $k$ \hspace{-.09in} & 200& 300 & 500 & 50 & 200 & 80 & 200& 500 & 4 \\
   \Xhline{1pt}

  \hspace{-.09in}   $t_1$ \hspace{-.09in}
       & \hspace{-.05in}    \textbf{32.1367} \hspace{-.09in}
       &  \hspace{-.05in}  \textbf{\textcolor{blue}{48.4165}} \hspace{-.09in}
       &  \hspace{-.05in}  79.8354 \hspace{-.09in}
        &  \hspace{-.05in}   10.5032 \hspace{-.09in}
         & \hspace{-.05in}   \textbf{26.3473} \hspace{-.09in}
         &  \hspace{-.05in}  \textbf{6.2474} \hspace{-.09in}
         &  \hspace{-.05in}  15.2432 \hspace{-.09in}
         &  \hspace{-.05in}  \textbf{147.4556} \hspace{-.09in}
         & \hspace{-.07in}     \textbf{0.3145}  \hspace{-.09in}    \\
               \hline

   \hspace{-.09in}  $t_2$  \hspace{-.09in}
       & \hspace{-.05in}   36.9786 \hspace{-.09in}
       &  \hspace{-.05in}  74.8509 \hspace{-.09in}
       &  \hspace{-.05in}  80.7513 \hspace{-.09in}
        & \hspace{-.05in}    \textbf{ 8.8483} \hspace{-.09in}
        & \hspace{-.05in}   29.7811 \hspace{-.09in}
        & \hspace{-.05in}   6.4447 \hspace{-.09in}
         & \hspace{-.05in}   15.2507 \hspace{-.09in}
         &  \hspace{-.05in}  155.4476 \hspace{-.09in}
         &  \hspace{-.07in}    0.3407 \hspace{-.09in}  \\
               \hline

    \hspace{-.09in}  $t_3$  \hspace{-.09in}
       &  \hspace{-.05in}  80.1901 \hspace{-.09in}
       &   \hspace{-.05in}       141.7107 \hspace{-.09in}
       & \hspace{-.05in}   \textbf{\textcolor{blue}{77.0104}} \hspace{-.09in}
       &  \hspace{-.05in}   65.2870 \hspace{-.09in}
       &  \hspace{-.05in}  63.9953 \hspace{-.09in}
       &  \hspace{-.05in}     11.5692 \hspace{-.09in}
       & \hspace{-.05in}   \textbf{\textcolor{blue}{11.7546}} \hspace{-.09in}
       &  \hspace{-.05in}  374.6021 \hspace{-.09in}
        &  \hspace{-.07in}  22.6884 \hspace{-.09in}  \\
                \Xhline{0.9pt}

 \hspace{-.09in}  $t_4(Z_1)$  \hspace{-.09in}
       & \hspace{-.05in}    \textbf{\textcolor{blue}{31.9943}} \hspace{-.09in}
       &    \hspace{-.05in}      \textbf{69.8041} \hspace{-.09in}
       &  \hspace{-.05in}  \textbf{79.6817} \hspace{-.09in}
       &  \hspace{-.05in}      9.5098 \hspace{-.09in}
        &  \hspace{-.05in}     26.3567 \hspace{-.09in}
        &  \hspace{-.05in}   6.1913 \hspace{-.09in}
        & \hspace{-.05in}   \textbf{15.2210} \hspace{-.09in}
        &  \hspace{-.05in}   149.8696 \hspace{-.09in}
           &  \hspace{-.07in}  stop at j=3 \hspace{-.09in}  \\
               \hline

  \hspace{-.09in}  $t_5(Z_1)$  \hspace{-.09in}
       &  \hspace{-.05in}  32.0460 \hspace{-.09in}
       &  \hspace{-.05in}  77.7833 \hspace{-.09in}
       & \hspace{-.05in}   81.3535 \hspace{-.09in}
       & \hspace{-.05in}    \textbf{9.2100} \hspace{-.09in}
       &  \hspace{-.05in}  \textbf{26.3251} \hspace{-.09in}
       & \hspace{-.05in}   \textbf{\textcolor{blue}{6.0000}} \hspace{-.09in}
        & \hspace{-.05in}     15.3378 \hspace{-.09in}
        &  \hspace{-.05in}  \textbf{\textcolor{blue}{144.2378}} \hspace{-.09in}
    &   \hspace{-.07in}  no\,Th.\!3.2(ii) \hspace{-.09in} \\
                           \Xhline{0.9pt}

 \hspace{-.09in}  $t_4(Z_2)$  \hspace{-.09in}
       & \hspace{-.05in}    36.4463 \hspace{-.09in}
       &  \hspace{-.05in}  70.1102 \hspace{-.09in}
       & \hspace{-.05in}   80.9745 \hspace{-.09in}
       &  \hspace{-.05in}   11.0663 \hspace{-.09in}
	    & \hspace{-.05in}    \textbf{\textcolor{blue}{26.2390}} \hspace{-.09in}
       & \hspace{-.05in}   7.0013 \hspace{-.09in}
        & \hspace{-.05in}   \textbf{15.3719} \hspace{-.09in}
        & \hspace{-.05in}   151.9858 \hspace{-.09in}
        & \hspace{-.07in}    0.2891 \hspace{-.09in}  \\
               \hline

 \hspace{-.09in}  $t_5(Z_2)$  \hspace{-.09in}
       & \hspace{-.05in}    \textbf{32.1301}  \hspace{-.09in}
       & \hspace{-.05in}   \textbf{57.7437} \hspace{-.09in}
       & \hspace{-.05in}   \textbf{ 79.2500} \hspace{-.09in}
       &  \hspace{-.05in}   \textbf{\textcolor{blue}{8.6892}} \hspace{-.09in}
        &  \hspace{-.05in}  30.0706 \hspace{-.09in}
        &  \hspace{-.05in}  \textbf{6.0525} \hspace{-.09in}
        &  \hspace{-.05in}  17.9037 \hspace{-.09in}
        &  \hspace{-.05in}  \textbf{150.2319} \hspace{-.09in}
        &   \hspace{-.07in}  \textbf{\textcolor{blue}{0.2522}} \hspace{-.09in}  \\
       \Xhline{0.9pt}

   \hspace{-.09in}   $e_1$  \hspace{-.09in}
       & \hspace{-.05in}   6.3344e-14 \hspace{-.09in}
       & \hspace{-.05in}   7.7079e-14 \hspace{-.09in}
       &  \hspace{-.05in}   1.1377e-13 \hspace{-.09in}
       &  \hspace{-.05in}  3.0507e-14 \hspace{-.09in}
       &  \hspace{-.05in}  2.1044e-13 \hspace{-.09in}
       & \hspace{-.05in}   1.5939e-14 \hspace{-.09in}
       &  \hspace{-.05in}  2.0370e-14 \hspace{-.09in}
       &  \hspace{-.05in}  2.1916e-13 \hspace{-.09in}
       &  \hspace{-.07in}   2.2204e-16 \hspace{-.09in} \\
               \hline

  \hspace{-.09in}   $e_2$  \hspace{-.09in}
       & \hspace{-.05in}   6.2526e-14 \hspace{-.09in}
       &  \hspace{-.05in}  7.6374e-14 \hspace{-.09in}
       &  \hspace{-.05in}   1.1156e-13 \hspace{-.09in}
       &  \hspace{-.05in}  3.0589e-14 \hspace{-.09in}
       &  \hspace{-.05in}    2.0469e-13 \hspace{-.09in}
        &  \hspace{-.05in}    1.5934e-14 \hspace{-.09in}
        &  \hspace{-.05in}  2.0358e-14 \hspace{-.09in}
        &  \hspace{-.05in}     2.1441e-13 \hspace{-.09in}
         &  \hspace{-.07in}    2.2204e-16 \hspace{-.09in}  \\
               \hline

   \hspace{-.09in}  $e_3$  \hspace{-.09in}
       &  \hspace{-.05in}    2.7159e-14 \hspace{-.09in}
       & \hspace{-.05in}   3.3145e-14 \hspace{-.09in}
        &   \hspace{-.05in}   3.2635e-14 \hspace{-.09in}
        &   \hspace{-.05in}    2.9561e-14 \hspace{-.09in}
        &  \hspace{-.05in}   4.1709e-14 \hspace{-.09in}
        & \hspace{-.05in}    1.4724e-14 \hspace{-.09in}
        &  \hspace{-.05in}  1.4372e-14 \hspace{-.09in}
        & \hspace{-.05in}   4.5437e-14 \hspace{-.09in}
        & \hspace{-.07in}   2.2204e-16 \hspace{-.09in} \\
               \hline

 \hspace{-.09in}  $e_4(Z_1)$  \hspace{-.09in}
       &  \hspace{-.05in}   2.8838e-04 \hspace{-.09in}
       &  \hspace{-.05in}  0.0010 \hspace{-.09in}
       &  \hspace{-.05in}    8.2355e-04 \hspace{-.09in}
       &  \hspace{-.05in}   1.5013e-05 \hspace{-.09in}
       & \hspace{-.05in}   2.5971e-04 \hspace{-.09in}
       & \hspace{-.05in}   4.5256e-06 \hspace{-.09in}
       &    \hspace{-.05in}  1.4281e-04 \hspace{-.09in}
        &   \hspace{-.05in}  2.3799e-04 \hspace{-.09in}
        &  \hspace{-.07in}   no \hspace{-.05in}  \\
               \hline

  \hspace{-.09in}  $e_5(Z_1)$  \hspace{-.09in}
       &   \hspace{-.05in}  6.8226e-06 \hspace{-.09in}
       &  \hspace{-.05in}  2.4236e-05 \hspace{-.09in}
       &  \hspace{-.05in}  4.8157e-06 \hspace{-.09in}
       &   \hspace{-.05in}   1.3170e-06 \hspace{-.09in}
       &  \hspace{-.05in}   8.7436e-07 \hspace{-.09in}
       &  \hspace{-.05in}  5.6164e-07 \hspace{-.09in}
        & \hspace{-.05in}   1.7502e-06 \hspace{-.09in}
        & \hspace{-.05in}   9.6339e-06 \hspace{-.09in}
        & \hspace{-.07in}   no \hspace{-.05in}  \\
               \hline

 \hspace{-.09in}  $e_4(Z_2)$  \hspace{-.09in}
       & \hspace{-.05in}   2.2986e-05 \hspace{-.09in}
       & \hspace{-.05in}  1.4722e-04 \hspace{-.09in}
        &  \hspace{-.05in}    0.0016 \hspace{-.09in}
        &  \hspace{-.05in}   2.0118e-06 \hspace{-.09in}
         &  \hspace{-.05in}   3.0671e-04 \hspace{-.09in}
         &   \hspace{-.05in}  4.0523e-05 \hspace{-.09in}
          &  \hspace{-.05in}  3.3912e-04 \hspace{-.09in}
          &   \hspace{-.05in}  0.0306 \hspace{-.09in}
           &  \hspace{-.07in}   4.4474e-13 \hspace{-.09in} \\
               \hline

  \hspace{-.09in}  $e_5(Z_2)$  \hspace{-.09in}
       &  \hspace{-.05in}    1.7899e-06 \hspace{-.09in}
       & \hspace{-.05in}   3.2424e-06 \hspace{-.09in}
       &  \hspace{-.05in}    1.9454e-05 \hspace{-.09in}
       &  \hspace{-.05in}    6.4818e-07 \hspace{-.09in}
       &   \hspace{-.05in}   7.2992e-06 \hspace{-.09in}
       & \hspace{-.05in}    4.8560e-07 \hspace{-.09in}
       &  \hspace{-.05in}  1.0330e-05 \hspace{-.09in}
       & \hspace{-.05in}  3.9119e-06 \hspace{-.09in}
       & \hspace{-.07in}  5.6066e-13 \hspace{-.09in} \\
               \hline

    \end{tabular}
    }
        \caption{Tests}
   \end{center}
\end{table}

The tests reflected in Table 1 show the following:
in certain cases, the iterative algorithm in  Theorem~\ref{thm2.5}(iii) outperforms Theorem \ref{thm2.5}(ii) (Tests 1-3,5-9)
and the converse also holds (Test 4); in certain cases, both of these approaches perform  (significantly) faster then
the pseudo-inverse approach using the MATLAB pinv-function  (Tests 1,2,4-6,8,9)  and hence also faster than the LTFAT-code for general frames;
in certain cases, some non-canonical duals provide a faster procedure in comparison to the use of the canonical dual, i.e., in certain cases Theorem~\ref{thm2.5_arb}(i) or Theorem~\ref{thm2.5_arb}(ii) performs faster than Theorem~\ref{thm2.5}(iii) or Theorem~\ref{thm2.5}(ii)  (Tests 1,4-6,8,9);
the execution time of the algorithms depends also much on the values of $k$, $N$, $r$, $N/r$.

In conclusion, we may say that when we have a dual frame for $(x_n)_{n=1}^N,$ the canonical one or any other,
the algorithms presented in this paper can be efficient for computing a dual frame for the reduced frame $(x_n)_{n\in E^C}$ and can be used to enrich and improve LTFAT. In certain cases, the iterative algorithm  outperforms the procedure involving a matrix inversion, which justifies its consideration.
In addition to this fact, let us also note that
 the iterative algorithm provides simultaneously a dual frame for all the erasure sets $E_j=\{1,2,\ldots,j\}$, $j=1,2,\ldots,k$, which gives flexibility for simultaneous use of multiple erasure sets. 
In certain cases, the use of a non-canonical dual frame outperforms the use of the canonical dual,
which also justifies the interest to non-canonical duals. 
The
size of the frame,  its redundancy, the dimension of the spaces, and the cardinality of the erasure set,
may have significant influence on the execution times of the considered procedures.
Further tests and deeper investigation of appropriate dual frames and
optimal values of the aforementioned parameters for efficient computational purposes of each method will be the task of a future work.


\section*{Acknowledgement}
The work on this topic was partly supported by the Scientific and Technological Cooperation project Austria--Croatia ``Frames, Reconstruction, and Applications'' (HR 03/2020).
The second author acknowledges support from the Austrian Science Fund (FWF) under grants P 35846-N and P 32055-N31.

\end{document}